\documentstyle[amscd,amssymb,verbatim,12pt]{amsart}
\theoremstyle{plain}
\newtheorem{Thm}{Theorem}[section]

\newtheorem{Prop}[Thm]{Proposition}
\theoremstyle{definition}
\newtheorem{Rem}[Thm]{Remark}
\newtheorem{Res}[Thm]{Result}
\newtheorem{Conj}[Thm]{Conjecture}

\numberwithin{equation}{section}
\newcommand{\eqnum}{\setcounter{equation}{\value{Thm}}}
\newcommand{\thnum}{\setcounter{Thm}{\value{equation}}}

\begin{document}
\title[BORDISM BETWEEN DOLD AND MILNOR MANIFOLDS]%
	{BORDISM BETWEEN DOLD AND MILNOR MANIFOLDS}
\author[A. K. DAS]%
	{ ASHISH KUMAR DAS}
\address{Department of Mathematics, North Eastern Hill University, Permanent 
Campus, Shillong-793022, Meghalaya, India.}  \email{akdas@@nehu.ac.in}
\thanks{}
\date{}
\begin{abstract} It is well known that Dold and Milnor manifolds give
generators for the unoriented bordism algebra ${\frak{N}}_*$ over ${\Bbb{Z}}
_2$. The purpose of this paper is to determine those Milnor manifolds which
represent the same bordism classes in ${\frak{N}}_*$ as their Dold
counterparts.  
\end{abstract} 
\maketitle

\section{Introduction} \label{S:intro}

Dold manifolds were first introduced by Dold \cite{aD01}. A Dold manifold $P(m,n)$ of dimension $m+2n$ is the quotient $({\Bbb{S}}^m \times {\Bbb{C}}P^n)/\sim$, where $\sim$ is an equivalence relation given by $(x,[z]) \sim (-x,[\bar{z}])$.  The ring structure of $H^{*}(P(m,n);{\Bbb{Z}} _2)$ is described as
\[
H^{*}(P(m,n);{\Bbb{Z}} _2) = \left[{\frac{{\Bbb{Z}} _2[c]}{c^{m+1} = 0}}\right] \otimes \left[{\frac{{\Bbb{Z}} _2[d]}{d^{n+1} = 0}}\right],
\]
and the total Stiefel-Whitney class  of $P(m,n)$ is given by 
\[
W(P(m,n)) = (1+c)^m (1+c+d)^{n+1},
\]
\noindent  where $c$ is the nonzero element of $H^{1}(P(m,n);{\Bbb{Z}} _2) \cong {\Bbb{Z}} _2$ and $d$ is a suitable nonzero element of $H^{1}(P(m,n);{\Bbb{Z}} _2) \cong {\Bbb{Z}} _2 \oplus {\Bbb{Z}} _2$.

On the other hand Milnor manifolds were first introduced by Milnor \cite{jM65}. A Milnor manifold $H(m,n)$ is the $m+n-1$ dimensional submanifold of ${\Bbb{R}}P^m \times {\Bbb{R}}P^n$ given by
\[
H(m,n) = \{([x_0, \dots ,x_m],[y_0, \dots ,y_n]) \in {\Bbb{R}}P^m \times {\Bbb{R}}P^n :\hspace{-.25cm} \sum_{i=0}^{\text {min} (m,n)} \hspace{-.25cm} x_i y_i = 0\}.
\]
\noindent In fact $H(m,n)$ is the submanifold of ${\Bbb{R}}P^m \times {\Bbb{R}}P^n$ dual to $(a+b)$;  $a$ and $b$ being the generators of $H^*({\Bbb{R}}P^m;{\Bbb{Z}} _2)$ and $H^*({\Bbb{R}}P^n;{\Bbb{Z}} _2)$ respectively. Note that $a^{m+1} = b^{n+1} = 0$, whereas $a^i$ and $b^j$ are non-zero for $0 \le i \le m$ and $0 \le j \le n$. The total Stiefel-Whitney class $W(H(m,n))$ of $H(m,n)$ is given by the restriction to $H(m,n)$ of the expression

\eqnum
\begin{equation} \label{E:sw}
{\frac{(1+a)^{m+1} (1+b)^{n+1}}{(1+a+b)}}
\end{equation}
\thnum

\noindent and the Stiefel-Whitney number of $H(m,n)$ corresponding to a partition $i_1 + i_2 + \dots + i_k = m+n-1$ is given by
\[
< (a+b)W_{i_1} \dots W_{i_k}, [{\Bbb{R}}P^m \times {\Bbb{R}}P^n] > \; \in {\Bbb{Z}} _2
\]
\noindent where $W_{i_j}$ is the sum of $i_j$ -dimensional terms in the expansion of the expression (\ref{E:sw}). Thus the partition $i_1 + i_2 + \dots + i_k = m+n-1$ corresponds to a non-zero Stiefel-whitney number of  $H(m,n)$ if and only if $W_{i_1} W_{i_2} \dots W_{i_k} = a^m b^{n-1}$ or $a^{m-1} b^n$. Since $H(m,n) \cong H(n,m)$, we always assume that $m \le n$. Further, throughout this paper the symbol ${ \hbox{ } \choose \hbox{ } }$  will stand for binomial coefficient reduced modulo $2.$ 

It is well known (see \cite{aD01},\cite{jM65}) that Dold as well as Milnor manifolds independently give generators for the unoriented bordism algebra ${\frak{N}}_*$ over ${\Bbb{Z}} _2$. Thus, in principle,  every Milnor manifold is bordant to a union of products of Dold manifolds. In this paper we investigate the question ``which Milnor manifolds are bordant to Dold manifolds?'' and obtain definite answers in most of the cases.

\section{Milnor manifolds which are bordant to Dold manifolds} \label{S:bor}

We begin with the following trivial remark

 \begin{Rem}\label{Rm:r1}
 $H(0,n) \cong {\Bbb{R}}P^{n-1} \cong P(n-1,0) \; \forall \; n \ge 1$.
\end{Rem}

For non-trivial situations we have 

\begin{Prop}\label{P:p1}
$H(2^{\alpha} - 2,n)$ is bordant to $P(n-2^{\alpha} +1, 2^{\alpha} -2) \; \; \forall \; \alpha $ $\ge 1$ \; and \; $\forall \; n \ge 2^{\alpha} -1$.
\end{Prop}

\begin{pf}
The total Stiefel-Whitney class of $P(n-2^{\alpha} +1, 2^{\alpha} -2)$ is given by 
\[
W(P) \; = \; (1+c)^{n-2^{\alpha} +1} (1+c+d)^{2^{\alpha} -1},
\]
\noindent where $c^{n-2^{\alpha} +2} \; = \; 0 \; = \; d^{2^{\alpha} -1}.$

On the other hand the total Stiefel-Whitney class of $H(2^{\alpha} -2,n)$ is given by

\begin{align*}
W(H) =&\;  {\frac{(1+a)^{2^{\alpha} -1} (1+b)^{n+1}}{(1+a+b)}}, \hspace{1cm} {\text{where}} \; a^{2^{\alpha} -1} = 0 = b^{n+1},\\
=&\;  (1+b)^{n-2^{\alpha} +1} (1+a)^{2^{\alpha} - 1} (1+b)^{2^{\alpha}} (1+a+b)^{2^{\alpha} - 1} (1+b^{2^{\alpha}}) \dots \\
&\;  (1+b^{2^{\gamma}}), \hspace{1cm} {\text{where}}\;  2^{\gamma} \le n < 2^{{\gamma} +1},\\
=&\;  (1+b)^{n-2^{\alpha} +1} (1 + b + a^{2} + ab)^{2^{\alpha} -1}.
\end{align*}

Consider a polynomial
\[
W(x,y) = (1+x)^{n-2^{\alpha} +1} (1+x+y)^{2^{\alpha} -1} \; \in {\Bbb{Z}} _2 [x,y]
\]
\noindent where $(x,y)$ satisfies the following conditions:

\eqnum
\begin{equation}\label{E:cn0}
{\text {dim} }\; x = 1, \; \; {\text {dim} }\; y = 2 , \; \; {\text {and} }\; \; y^{2^{\alpha} -1} = 0.
\end{equation}
\thnum

\noindent Then, 
\[
W(P) = W(c,d) \; {\text{and}} \; W(H) = W(b,a(a+b)).
\]
\noindent Clearly, both  $(c,d)$ and $(b,a(a+b))$ satisfy the conditions in (\ref{E:cn0}). Let \; 
$W_i (x,y)$ be the sum of $i$-dimensional terms in $W(x,y)$, where $i > 0$. Then
\[
W_i (x,y) = D_{0,i} x^i + D_{1,i} x^{i-2} y + \dots + D_{2^{\alpha} -2,i} x^{i-2^{\alpha +1} + 4} y^{2^{\alpha} -2}
\]
\noindent where $D_{j,i} \in {\Bbb{Z}} _2$ ,   $0 \le j \le 2^{\alpha} -2$, \; and \; $D_{j,i} = 0$ \; if  $i < 2j$.  Let $ \omega \equiv i_1 + i_2 + \dots + i_k \; = \; n+ 2^{\alpha} -3$ be a partition of $n+ 2^{\alpha} -3$. Then,
\[
\textstyle
W_{\omega} (x,y) = W_{i_{1}} (x,y)  \dots W_{i_{k}} (x,y)  =  \sum_{t=0}^{2^{\alpha} -2}{f_t (D) x^{n+2^{\alpha} -3 -2t} y^t}
\]
\noindent where
\[
f_t (D) = \sum_{j_i + \dots + j_k = t}{D_{j_1,i_1}D_{j_2,i_2} \dots D_{j_k,i_k}} \; \in  {\Bbb{Z}} _2.
\]
\noindent Therefore,
\[
W_{\omega} (c,d) = f_{2^{\alpha} -2} (D) c^{n-2^{\alpha} +1} d^{2^{\alpha} -2}.
\]

  On the other hand
\[
W_{\omega} (b,a(a+b)) = \sum_{t=0}^{2^{\alpha} -2}{f_t (D) b^{n+2^{\alpha} -3 -2t} a^t (a+b)^t}.
\]
\noindent Let $X_t = b^{n+2^{\alpha} -3 -2t} a^t (a+b)^t$ \; where $0 \le t \le 2^{\alpha} -2.$ \; If $ t = 2^{\alpha} -2$, then 
\begin{equation*}
X_t = \; b^{n-2^{\alpha} +1} a^{2^{\alpha} -2} (a+b)^{2^{\alpha} -2}
= \; a^{2^{\alpha} -2} b^{n - 1}.
\end{equation*}

\noindent If $ t \le 2^{\alpha -1} -2$, then $n+2^{\alpha} -3 -2t \ge n+2^{\alpha} -3 -(2^{\alpha} -4) = n+1,$ and so \; $X_t = 0$ \; (since $b^{n+1} = 0$). \; So, let $2^{\alpha -1} -1 \le t < 2^{\alpha} -2.$  Then,

\begin{align*}
X_t &= \;  {t \choose {2^{\alpha} -2-t}} a^{2^{\alpha} -2} b^{n-1} + {t \choose {2^{\alpha} -2 -t-1}} a^{2^{\alpha} -3} b^n\\
&= \; A_t (a^{2^{\alpha} -2} b^{n-1} + a^{2^{\alpha} -3} b^n), 
\end{align*}

\noindent where \; $A_t = 0$ or $1$; noting that
\[
{t \choose {2^{\alpha} -2-t}} + {t \choose {2^{\alpha} -2-t-1}} = {t+1 \choose {2^{\alpha} -2-t}} = 0 \in  {\Bbb{Z}} _2.
\]
\noindent Therefore, 
\[
W_{\omega} (b,a(a+b)) = f_{2^{\alpha} -2} (D) a^{2^{\alpha} -2} b^{n-1} + A (a^{2^{\alpha} -2} b^{n-1} + a^{2^{\alpha} -3} b^n)
\]
\noindent where $A = 0$ or $1$. Hence it follows that the stiefel-whitney number
\begin{align*}
&\; \langle W_{\omega} (P), [P] \rangle \; = \; \langle W_{\omega} (c,d), [P] \rangle
\; = \; f_{2^{\alpha} -2} (D)\\
=&\; \langle W_{\omega} (b,a(a+b)), [H] \rangle \;
= \; \langle W_{\omega} (H), [H] \rangle
\end{align*}

\noindent This completes the proof.
\end{pf}

\begin{Prop}\label{P:p2}
$H(m, m + 2^{\alpha} B)$ is bordant to $P(2^{\alpha} B -1,m) \; \; \forall \; m$ $ \ge 0, \; \forall \; B \ge 1, and \; \forall \; \alpha $ such that  $2^{\alpha} \ge m+1$.
\end{Prop}

\begin{pf}
The total Stiefel-Whitney class of $P(2^{\alpha} B -1, m)$ is given by 
\[
W(P) \; = \; (1+c)^{2^{\alpha} B -1} (1+c+d)^{m +1},
\]
\noindent where $c^{2^{\alpha} B} \; = \; 0 \; = \; d^{m +1}.$

On the other hand the total Stiefel-Whitney class of $H(m, m+ 2^{\alpha} B)$ is given by

\begin{align*}
W(H) =& \;  {\frac{(1+a)^{m +1} (1+b)^{m + 2^{\alpha } B +1}}{(1+a+b)}}, \quad \quad {\text{where}} \; a^{m +1} = 0 = b^{m +2^{\alpha} B +1},\\
=& \;  (1+a+b+ab)^{m + 1} (1+b)^{2^{\alpha} B} (1+a+b)^{2^{\alpha} -1} (1+b^{2^{\alpha}}) \dots\\
&\;  (1+b^{2^{\gamma}}), \hspace{1cm}  {\text{where}}\;  2^{\gamma} \le m+ 2^{\alpha} B < 2^{{\gamma} +1},\\
=&\; (1+a+b+ab)^{m + 1} (1+a+b)^{2^{\alpha} -1} (1+ b)^{2^{\alpha} (B-1)}\\
=&\; (1+a+b+ab)^{m + 1} (1+a+b)^{2^{\alpha} B -1}, \hspace{1cm} {\text{since $a^{2^{\alpha}} = 0$.}}
\end{align*}

Let $W(x,y)$ $=   (1+x)^{2^{\alpha} B -1} (1+x+y)^{m + 1} \; \in {\Bbb{Z}} _2 [x,y]$, \; where $(x,y)$ satisfies the following conditions:

\eqnum
\begin{equation}\label{E:cn1}
{\text {dim} }\; x = 1, \;  {\text {dim} }\; y = 2 , \;  {\text {and} }\;
\; x^{2m + 2^{\alpha} B}  =  0  =  y^{m + 1}.
\end{equation}
\thnum

\noindent Then 
\[
W(P) = W(c,d),\;  \; {\text{and}} \; \; W(H) = W(a+b, ab).
\]
\noindent Clearly, the pairs $(c,d)$ and $(a+b,ab)$  satisfy the conditions in (\ref{E:cn1}).

Let \; $W_i (x,y)$ be the sum of $i$-dimensional terms in $W(x,y)$, where $i > 0$. Then,
\[
W_i (x,y) = D_{0,i} x^i + D_{1,i} x^{i-2} y + \dots + D_{m,i} x^{i-2m} y^m ,
\]
\noindent where $D_{j,i}  \;  \in \; {\Bbb{Z}} _2 \; \forall \; j = 0,1, \dots ,m.$   Let $ \omega \equiv i_1 + i_2 + \dots + i_k \; = \; 2m + 2^{\alpha} B -1$ be a partition of $ 2m + 2^{\alpha} B -1$. Then, using conditions in (\ref{E:cn1}),
\begin{align*}
&\; W_{\omega} (x,y) = W_{i_{1}} (x,y)  \dots W_{i_{k}} (x,y) \\
=&\;  \sum_{t=0}^{m} {f_t (D_{0, i_1} \dots D_{m, i_k} ) x^{2m +2^{\alpha} B -1 -2t} y^t}
\end{align*}
\noindent where $f_t (D_{0, i_1} \dots D_{m, i_k}) \in {\Bbb{Z}} _2$  is a polynomial in  $D_{j,i_r}, \; 0 \le j \le m, \; 1 \le r \le k.$ \; Therefore,
\[
W_{\omega} (c,d) = f_m (D_{0, i_1} \dots D_{m, i_k}) c^{2^{\alpha} B -1} d^m, \; {\text{ since $c^{2^{\alpha} B} = 0$.}}
\]

On the other hand 
\[
W_{\omega} (a+b,ab) =  \sum_{t=0}^{m} {f_t (D_{0, i_1} \dots D_{m, i_k} ) (a+b)^{2m +2^{\alpha} B -1 -2t} (ab)^t}.
\]
\noindent Let $X_t = (a+b)^{2m +2^{\alpha} B -1 -2t} (ab)^t$ \; where $0 \le t \le m$.  If $t=m,$ then
\[
X_t =\; (a+b)^{2^{\alpha} B-1} (ab)^m =\; a^m b^{m + 2^{\alpha} B -1},
\]
\noindent and if $t <m,$ then
\begin{align*}
X_t =&\; (a+b)^{2m +2^{\alpha} B -1 -2t} (ab)^t\\
=&\; {2m +2^{\alpha} B -1 -2t \choose m-t} a^m b^{m +2^{\alpha} B -1} \; + \\
&\; {2m +2^{\alpha} B -1 -2t \choose m-t-1} a^{m-1} b^{m +2^{\alpha} B}\\
=&\; A_t (a^m b^{m +2^{\alpha} B -1} + a^{m-1} b^{m +2^{\alpha} B}),
\end{align*}
\noindent where $A_t = 0$ or $1$; \; noting that 
\begin{align*}
&\; {2m +2^{\alpha} B -1 -2t \choose m-t} \; + \; {2m +2^{\alpha} B -1 -2t \choose m-t-1} \\
=&\; {2m +2^{\alpha} B  -2t \choose m-t} \; = \; 0  \in  {\Bbb{Z}} _2, 
\end{align*}
\noindent looking at the lowest power of $2$ in $(m-t)$. Thus,
\begin{align*}
&\; W_{\omega} (a+b,ab)\\
=&\;  f_m (D_{0, i_1} \dots D_{m, i_k} ) a^m b^{m +2^{\alpha} B -1} +  A (a^m b^{m +2^{\alpha} B -1} + a^{m-1} b^{m +2^{\alpha} B})
\end{align*}
\noindent where $A = 0$ or $1$.  Hence, it follows that the Stiefel-Whitney number
\begin{align*}
&\; \langle W_{\omega} (P), [P] \rangle =  \langle W_{\omega} (c,d), [P] \rangle = f_m (D_{0, i_1} \dots D_{m, i_k} ) \\
=&\;  \langle W_{\omega} (a+b,ab), [H] \rangle =\;  \langle W_{\omega} (H), [H] \rangle.
\end{align*}
\noindent This completes the proof.
\end{pf}

\begin{Prop}\label{P:p3}
$H(2^{\alpha}, 2^{\alpha +1} B)$ is bordant to $P(2^{\alpha +1} B -2^{\alpha} - 1, 2^{\alpha})$
$\forall$$  \; \alpha \ge 1$ and \; $\forall \; B \ge 1$.
\end{Prop}
\begin{pf}
The total Stiefel-Whitney class of $P(2^{\alpha +1} B- 2^{\alpha} -1, 2^{\alpha})$ is given by 
\[
W(P) \; = \; (1+c)^{2^{\alpha +1} B- 2^{\alpha} -1} (1+c+d)^{2^{\alpha} +1},
\]
\noindent where $c^{2^{\alpha +1} B- 2^{\alpha} } \; = \; 0 \; = \; d^{2^{\alpha} +1}.$

On the other hand the total Stiefel-Whitney class of $H(2^{\alpha} ,2^{\alpha +1} B)$ is given by
\begin{align*}
W(H) =& \;  {\frac{(1+a)^{2^{\alpha} +1} (1+b)^{2^{\alpha +1} B +1}}{(1+a+b)}}, \hspace{1cm} {\text{where}} \; a^{2^{\alpha} +1} = 0 = b^{2^{\alpha +1} B +1},\\
=& \;  (1+a)^{2^{\alpha} + 1} (1+b)^{2^{\alpha +1} B +1} (1+a+b) \dots (1+a^{2^{\alpha}} + b^{2^{\alpha}})\\
&\; (1+b^{2^{\alpha + 1}}) \dots  (1+b^{2^{\gamma}}), \hspace{1cm}  {\text{where}}\;  2^{\gamma} \le 2^{\alpha +1} B < 2^{{\gamma} +1},\\
=&\; (1+a)^{2^{\alpha} + 1} (1+b)^{2^{\alpha +1} (B-1)+1} (1+a+b)^{2^{\alpha} -1} (1+a^{2^{\alpha}} + b^{2^{\alpha}})\\
=&\; (1+a)^{2^{\alpha} + 1} (1+b)^{2^{\alpha +1} (B-1)+1} (1+a+b)^{2^{\alpha} -1} (1 + b^{2^{\alpha}})\\
&\; + (1+a)^{2^{\alpha} + 1} (1+b)^{2^{\alpha +1} (B-1)+1} (1+a+b)^{2^{\alpha} -1} a^{2^{\alpha}}\\
=&\;  (1+a+b+ab)^{2^{\alpha} + 1} (1+a+b)^{2^{\alpha +1} (B-1)} (1+a+b)^{2^{\alpha} -1} \\
&\; + (1+b)^{2^{\alpha +1} (B-1)+1} (1+b)^{2^{\alpha} -1} a^{2^{\alpha}}, \quad \; {\text{since}}\; a^{2^{\alpha}+1} = 0,\\
=&\; (1+a+b+ab)^{2^{\alpha} + 1} (1+a+b)^{2^{\alpha +1} B - 2^{\alpha} -1} \\
&\; + (1+b)^{2^{\alpha +1} B - 2^{\alpha}} a^{2^{\alpha}}.
\end{align*}
Let $W(x,y)$ $= \;  (1+x)^{2^{\alpha +1} B - 2^{\alpha} -1} (1+x+y)^{2^{\alpha} + 1} \in {\Bbb{Z}} _2 [x,y]$, \; where $(x,y)$ satisfies the following conditions:

\eqnum
\begin{equation}\label{E:cn2}
{\text {dim} }\; x = 1, \;  {\text {dim} }\; y = 2 , \;  {\text {and} }\;
\; x^{2^{\alpha +1} B + 2^{\alpha}}  =  0  =  y^{2^{\alpha} + 1}.
\end{equation}
\thnum

\noindent Then 
\begin{align*}
W(P) &= W(c,d) {\text{\; and}}\\
W(H) &= W(a+b, ab) + a^{2^{\alpha}} (1+b)^{2^{\alpha +1} B - 2^{\alpha}}.
\end{align*}
\noindent Clearly, the pairs $(a+b,ab)$ and $(c,d)$ satisfy the conditions in (\ref{E:cn2}), noting that the term $a^{2^{\alpha}} b^{2^{\alpha +1} B}$ may be omitted due to dimensional considerations.

Let \; $W_i (x,y)$ be the sum of $i$-dimensional terms in $W(x,y)$, where $i > 0$. Now,
\[
W(x,y) =\; (1+x)^{2^{\alpha +1} B } + (1+x)^{2^{\alpha +1} B - 1} y + (1+x)^{2^{\alpha +1} B - 2^{\alpha} } y^{2^{\alpha}}.
\]
\noindent Therefore,
\[
W_i (x,y) = D_i x^i + E_i x^{i-2} y + F_i x^{i-2^{\alpha +1}} y^{2^{\alpha}}
\]
\noindent where
\[
D_i = \; {2^{\alpha +1} B \choose i}, \; E_i = \; {2^{\alpha +1} B - 1 \choose i-2}, \; F_i =\;{2^{\alpha +1} B - 2^{\alpha} \choose i- 2^{\alpha +1}}, 
\]
\noindent using the convention that ${r \choose s} = 0$ if $s<0$ or $s>r$.  Let $\delta$ be the exponent of $2$ in $2^{\alpha + 1} B$, $i.e.$ \; $2^{\alpha + 1} B$ $= \; 2^{\delta}.$odd.  Clearly, $ \delta \ge \alpha +1.$  Now,
\[
i(i-1)D_i = (2^{\alpha + 1} B +1 -i) 2^{\alpha + 1} B E_i.
\]
\noindent So, it follows that if $i =\; 2^{\beta}.$odd, with $\beta \ge 1$, then
\begin{align*}
\beta = \delta \Rightarrow &D_i = E-i,\\
\beta < \delta \Rightarrow &D_i = 0,\\
\beta > \delta \Rightarrow &E-i = 0.
\end{align*}
\noindent Further note that $D_i = 0 = F_i$ \; if $i$ is odd.  Therefore,
\[
W_i(x,y) =
\begin{cases}
E_i x^{i-2} y, \; \; \; \text{if $i$ is odd,}\\
E_i x^{i-2} y + F_i x^{i-2^{\alpha +1}} y^{2^{\alpha}}, \; \; \; \text{if $i=2^{\beta}.$odd, $1 \le \beta < \delta$,}\\
D_i (x^i + x^{i-2} y) + F_i x^{i-2^{\alpha +1}} y^{2^{\alpha}}, \; \; \; \text{if $i=2^{\delta}.$odd,}\\
D_i x^i + F_i x^{i-2^{\alpha +1}} y^{2^{\alpha}}, \; \; \; \text{if $i=2^{\beta}.$odd, $\beta > \delta.$}
\end{cases}
\]
\noindent  Let $ \omega \equiv i_1 + i_2 + \dots + i_k \; = \; 2^{\alpha + 1} B + 2^{\alpha} -1$ be a partition of $2^{\alpha + 1} B + 2^{\alpha} -1$. Clearly, at least one $i_j$ must be odd.   So, at least one $W_{i_j} (x,y)$ is a multiple of $y$.  Let $i_1, i_2, \dots ,i_m$ be of type $2^{\beta}.$odd, with $0 \le \beta < \delta.$  Then, $m \ge 1.$  Let $i_{m+1}, i_{m+2}, \dots , i_{m+n}$ be of type $2^{\delta}.$odd, and let $i_{m+n+1}, i_{m+n+2}, \dots , i_{k}$ be of type $2^{\beta}.$odd, with $\beta > \delta.$ Then, we have
\begin{align*}
&\; W_{\omega} (x,y) = W_{i_1} (x,y) \dots W_{i_k} (x,y)\\
=&\; \left( \underset{j=1}{\overset{m}{\prod}} E_{i_j} x^{i_j -2} y \right) \left( \underset{j=m+1}{\overset{m+n}{\prod}} D_{i_j} (x^{i_j} + x^{i_j -2} y) \right) \left( \underset{j=m+n+1}{\overset{k}{\prod}} D_{i_j} x^{i_j} \right)\\
=&\;  \underset{t=0}{\overset{n}{\sum}} D.E.G_t x^{2^{\alpha + 1} B + 2^{\alpha} -1 -2(m+t)} y^{m+t},
\end{align*}
\noindent where $D = \underset{j=m+1}{\overset{k}{\prod}} D_{i_j},$ \; $E = \underset{j=1}{\overset{m}{\prod}} E_{i_j},$ \; $G_t = {n \choose t}$, \; $0 \le t \le n.$  Thus,
\[
W_{\omega}(c,d) =
\begin{cases}
D.E.G_{2^{\alpha} - m} c^{2^{\alpha +1} B -2^{\alpha} -1} d^{2^{\alpha}}, \; \; \; \text{if $0 \le 2^{\alpha} -m \le n$}\\
0, \; \; \; {\text{otherwise}}.
\end{cases}
\]

 On the other hand
\[
W_{\omega} (a+b,ab) =  \underset{t=0}{\overset{n}{\sum}} D.E.G_t (a+b)^{2^{\alpha + 1} B + 2^{\alpha} -1 -2(m+t)} (ab)^{m+t}.
\]
\noindent Let $X_t = (a+b)^{2^{\alpha + 1} B + 2^{\alpha} -1 -2(m+t)} (ab)^{m+t},$ \; where $0 \le t \le n.$  Clearly, if $m+t > 2^{\alpha}$ then $X_t = 0$. If $m+t = 2^{\alpha}$ then
\[
X_t =\; (a+b)^{2^{\alpha + 1} B + 2^{\alpha} -1} (ab)^{2^{\alpha}} =\; a^{2^{\alpha}} b^{2^{\alpha + 1} B -1}.
\]
\noindent If $m+t < 2^{\alpha}$ then
\begin{align*}
X_t =&\; {2^{\alpha + 1} B + 2^{\alpha} -1 -2(m+t) \choose 2^{\alpha} - (m+t)} a^{2^{\alpha}} b^{2^{\alpha + 1} B -1} + \\
&\; {2^{\alpha + 1} B + 2^{\alpha} -1 -2(m+t) \choose 2^{\alpha} - (m+t) -1} a^{2^{\alpha} -1} b^{2^{\alpha + 1} B }\\
=&\; A_t ( a^{2^{\alpha}} b^{2^{\alpha + 1} B -1} + a^{2^{\alpha} -1} b^{2^{\alpha + 1} B }),
\end{align*}
\noindent where $A_t = 0$ or $1$;  noting that 
\begin{align*}
&\; {2^{\alpha + 1} B + 2^{\alpha} -1 -2(m+t) \choose 2^{\alpha} - (m+t)}  +  {2^{\alpha + 1} B + 2^{\alpha} -1 -2(m+t) \choose 2^{\alpha} - (m+t) -1}\\
=&\; {2^{\alpha + 1} B + 2^{\alpha} -2(m+t) \choose 2^{\alpha} - (m+t)} =\; 0 \;  \in \; {\Bbb{Z}} _2 ,
\end{align*}
\noindent looking at the lowest power of $2$ in $(m+t).$  Hence it follows that
\begin{align*}
&\; W_{\omega}(a+b,ab) \\
=&\; \begin{cases}
D.E.(G_{2^{\alpha} - m} a^{2^{\alpha}} b^{2^{\alpha +1} B  -1} + A.g(a,b)),  \text{if $0 \le 2^{\alpha} -m \le n,$}\\
A'.g(a,b), \; \;  {\text{ otherwise,}}
\end{cases}
\end{align*}
\noindent where $g(a,b)$ $= a^{2^{\alpha}} b^{2^{\alpha +1} B  -1} + a^{2^{\alpha} -1} b^{2^{\alpha +1} B }$, \; and $A, A' \; \in {\Bbb{Z}} _2 .$  Thus, we have
\begin{align*}
\langle W_{\omega} (c,d), [P] \rangle =&\; \begin{cases} 
D.E.G_{2^{\alpha} - m}, \; \text{if $0 \le 2^{\alpha} -m \le n$}\\
0, \; \; {\text{otherwise}}
\end{cases}\\
=&\; \langle W_{\omega} (a+b,ab), [H] \rangle.
\end{align*}

Now,  $W_i (H) = W_i (a+b.ab) + N_i (a,b)$, \; where 
\[
N_i (a,b) = {2^{\alpha +1} B -2^{\alpha} \choose i -2^{\alpha}} a^{2^{\alpha}} b^{i -2^{\alpha}}.
\]
\noindent Note that
\begin{align*}
(i)&\; N_i (a,b) N_j (a,b) = 0,\; \forall \; i, j \ge 1,\\
(ii)&\; N_i (a,b) W_j(a+b,ab) = 0, \; \forall \; i \ge 1,\; {\text{if $j$ is odd, \; and}}\\
(iii)&\; N_i (a,b) = 0, \;  {\text{if $i$ is odd}}.
\end{align*}

Hence, it follows that
\begin{align*}
 W_{\omega} (H) =&\; W_{i_{1}} (H)  \dots W_{i_{k}} (H)\\
=&\; W_{\omega} (a+b,ab) \; + \; \underset{j=1}{\overset{k}{\sum}} N_{i_{j}} (a,b) W_{i_{1}} (a+b,ab) \dots  \\
&\;  \widehat{W_{i_{j}}} (a+b,ab) \dots W_{i_{k}} (a+b,ab)\\
=&\; W_{\omega} (a+b,ab), \; {\text{ since $\omega$ has an odd summand}.}
\end{align*}
\noindent Thus , the Stiefel-Whitney number
\[
 \langle W_{\omega} (P), [P] \rangle =  \langle W_{\omega} (c,d), [P] \rangle = \langle W_{\omega} (a+b,ab), [H] \rangle  =  \langle W_{\omega} (H), [H] \rangle.
\]
\noindent This completes the proof.
\end{pf}

In \cite{jM65} it has been proved that 

\begin{Res}\label{Rs:r1}
${\frak{N}}_*$ is a polynomial ring over ${\Bbb{Z}} _2$ with independent generators ${\Bbb{R}}P^{2t}$ and $H(2^k, 2t2^k)$ where $ t, k \ge 1$.
\end {Res}

Therefore, inview of Remark \ref{Rm:r1} and Proposition \ref{P:p3} we have

\begin{Rem}\label{Rm:r2}
The Milnor manifolds lying in the generating set of ${\frak{N}}_*$ over ${\Bbb{Z}} _2$ are all bordant to Dold manifolds.
\end{Rem}

The following proposition is not directly related to the question we are investigating. However it has its own significance.

\begin{Prop}\label{P:p4}
$H(m, 2^{\alpha} -1)$ is bordant to  ${\Bbb{R}}P^m \times {\Bbb{R}}P^{2^{\alpha} -2} \; \; \forall \; m \ge 0, \; $ and $\; \forall \;  \alpha $ such that  $2^{\alpha} > m+1$.
\end{Prop}

\begin{pf}
The total Stiefel-Whitney class of $H(m, 2^{\alpha} -1)$ is given by
\begin{align*}
W(H) =& \;  {\frac{(1+a)^{m +1} (1+b)^{2^{\alpha }}}{(1+a+b)}}, \hspace{1cm} {\text{where}} \; a^{m +1} = 0 = b^{2^{\alpha} },\\
=& \;  (1+a)^{m + 1} (1+a+b)^{2^{\alpha} -1}.
\end{align*}

On the other hand the total Stiefel-Whitney class of the product $M=$ ${\Bbb{R}}P^m \times {\Bbb{R}}P^{2^{\alpha} -2}$ is given by 
\[
W(M) = (1+u)^{m + 1} (1+v)^{2^{\alpha} -1}.
\]
\noindent where $u \in H^*({\Bbb{R}}P^m; {\Bbb{Z}} _2)$ , \; $v \in H^*({\Bbb{R}}P^{2^{\alpha} -2}; {\Bbb{Z}} _2)$ are the generators of the respective cohomology rings.  Let
\[
W(x,y) = (1+x)^{m + 1} (1+y)^{2^{\alpha} -1},
\]
\noindent where $(x,y)$ satisfies the following conditions:
\eqnum
\begin{equation}\label{E:cn3}
{\text {dim} }\; x =  \; {\text {dim} }\; y = \; 1 , \; \; {\text {and} }\; \; x^{m+1} = 0.
\end{equation}
\thnum
\noindent Then 
\[
 W(H) = W(a, a+b) \; {\text{and}} \; W(M) = W(u, v).
\]
\noindent Clearly, the pairs $(a,a+b)$ and $(u,v)$  satisfy the conditions in (\ref{E:cn3}).
Let \; $W_i (x,y)$ be the sum of $i$-dimensional terms in $W(x,y)$, where $i > 0$. Then, $W_i (x,y)$ is a homogeneous polynomial in $x, y$ of dimension $i$.  Let $ \omega \equiv i_1 + i_2 + \dots + i_k \; = \; m + 2^{\alpha}  -2$ be a partition of $ m + 2^{\alpha} -2$. Then, using conditions in (\ref{E:cn3}),
\[
 W_{\omega} (x,y) = W_{i_{1}} (x,y)  \dots W_{i_{k}} (x,y) =  \sum_{j=0}^{m} {D_j x^j y^{m + 2^{\alpha} -2 -j}} 
\]
\noindent where $D_j  \;  \in \; {\Bbb{Z}} _2 \; \forall \; j = 0,1, \dots ,m.$  Therefore,
\[
W_{\omega} (M) =\; W_{\omega} (u,v) =\; D_m u^m v^{2^{\alpha} -2}.
\]

On the other hand,
\[
 W_{\omega} (H) =\;W_{\omega} (a, a+b) =\;  \sum_{j=0}^{m} {D_j a^j (a+b)^{m + 2^{\alpha} -2 -j}}.
\] 
\noindent Let $X_t = a^j (a+b)^{m +2^{\alpha} -2 -j}$ \; where $0 \le j \le m$. Then \; 
$X_m = a^m b^{2^{\alpha} -2},$ \; and if $j < m$,
\begin{align*}
X_j =&\; {m +2^{\alpha} -2 -j \choose m-j} a^m b^{2^{\alpha} -2}  +  {m +2^{\alpha} -2 -j \choose m-j-1} a^{m-1} b^{2^{\alpha} -1}\\
=&\; A_j (a^m b^{2^{\alpha} -2} + a^{m-1} b^{2^{\alpha} -1}),
\end{align*}
\noindent where $A_j = 0$ or $1$; \; noting that 
\begin{align*}
&\; {m +2^{\alpha} -2 -j \choose m-j} \; + \; {m +2^{\alpha} -2 -j \choose m-j-1} \\
=&\; {2^{\alpha} -1 +m -j \choose m-j} \; = \; 0  \in  {\Bbb{Z}} _2,  
\end{align*}
\noindent looking at the lowest power of $2$ in $(m-j)$. Thus,
\[
W_{\omega} (H) = D_m a^m b^{2^{\alpha} -2} +  A (a^m b^{2^{\alpha} -2} + a^{m-1} b^{2^{\alpha} -1}),
\]
\noindent where $A = 0$ or $1$.  Hence, it follows that the Stiefel-Whitney number
\[
\langle W_{\omega} (H), [H] \rangle = D_m = \langle W_{\omega} (M), [{\Bbb{R}}P^m \times {\Bbb{R}}P^{2^{\alpha} -2} ] \rangle  .
\]
\noindent This completes the proof.
\end{pf}

\section{Milnor manifolds which are not bordant to Dold manifolds} \label{S:nbor}

First of all note that the Milnor manifolds, which are boundaries, are trivially bordant to Dold manifolds which are also boundaries (of course of the same dimensions) and in that context we have the following results from \cite{kD00} and \cite{sK89}

\begin{Res}\label{Rs:r2}
A Milnor manifold $H(m,n)$, with $m \leq
n$, bounds if and only if at least one of the following conditions holds:

\noindent (a) $m = n$,

\noindent (b) $m = 1$,

\noindent (c) $mn \equiv 1$ (mod $2$),

\noindent (d) $n \equiv 2$(mod $4$) and $m+1 < 2^{\nu (n+2)},$ where $\nu
(n+2)$ is the largest integer such that ${2^{\nu (n+2)}| {(n+2)}}.$
\end{Res}

\begin{Res}\label{Rs:r3}
A Dold manifold $P(m,n)$ bounds if and only if one of the following conditions holds:\\
\noindent (a) $n$ is odd,\\
\noindent (b) $n$ is even, $m$ is odd, $m>n$ and $2^{\nu (m-n-1)} > n$.
\end{Res}

The class of Milnor manifolds which are not bordant to Dolds manifolds is much bigger than the class of Milnor manifolds which are actually bordant to Dolds manifolds, and it is quite evident from the following proposition

\begin{Prop}\label{P:p5}
Let $m$ be odd, and $n$ be even such that $ m < n$. Then $H(m,n)$ is not bordant to a Dold manifold unless $H(m,n)$ itself is a boundary.
\end{Prop}

\begin{pf}
Assume that $H(m,n)$ is not a boundary. We know (see \cite{jM65}) that $H(m,n)$ fibres over ${\Bbb{R}}P^m$ with fibre  ${\Bbb{R}}P^{n-1}$.  Therefore, the mod $2$ Euler characteristic 
\[
\chi (H(m,n)) \; = \; \chi ({\Bbb{R}}P^m ) \chi ({\Bbb{R}}P^{n-1} ) \; = \; 0.
\]

Let $d =$ \; dim $(H(m,n)) \;  = \; m+n-1$. Clearly $d$ is even and so the non-bounding Dold manifolds of dimension $d$ are of the type $P(r,s)$ \; where both $r$ and $s$ are even (by Result \ref{Rs:r3}) and $r+2s =d$. Since (see \cite{cF65}) $P(r,s)$ fibres over ${\Bbb{R}}P^r$ with fibre ${\Bbb{C}}P^s$, it follows that the mod $2$ Euler characteristic   
\[
\chi (P(r,s)) \;  = \; \chi ({\Bbb{R}}P^r ) \chi ({\Bbb{C}}P^s ) \;  \ne \;  0.
\]
\noindent Thus, $H(m,n)$ is not bordant to a Dold manifold.
\end{pf}

At this moment we can not say much about the remaining Milnor manifolds; however the Results \ref{Rs:r2} and \ref{Rs:r3} together with some computer calculations suggest us to make the following conjecture:

\begin{Conj}
The non-bounding Milnor manifolds which are not considered in Propositions \ref{P:p1}, \ref{P:p2}, \ref{P:p3}, \ref{P:p4} and \ref{P:p5}, are not bordant to Dold manifolds.
\end{Conj}

\end{document}